# Even Faster Accelerated Coordinate Descent Using Non-Uniform Sampling

Zeyuan Allen-Zhu
zeyuan@csail.mit.edu
Princeton University

Zheng Qu
zhengqu@hku.hk
University of Hong Kong

Peter Richtárik
peter.richtarik@ed.ac.uk
University of Edinburgh

Yang Yuan
yangyuan@cs.cornell.edu
Cornell University

December 30, 2015[*]

**Abstract**

Accelerated coordinate descent is widely used in optimization due to its cheap per-iteration cost and scalability to large-scale problems. Up to a primal-dual transformation, it is also the same as accelerated stochastic gradient descent that is one of the central methods used in machine learning.

In this paper, we improve the best known running time of accelerated coordinate descent by a factor up to $\sqrt{n}$. Our improvement is based on a clean, novel non-uniform sampling that selects each coordinate with a probability proportional to the square root of its smoothness parameter. Our proof technique also deviates from the classical estimation sequence technique used in prior work. Our speed-up applies to important problems such as empirical risk minimization and solving linear systems, both in theory and in practice.

---

[*]The results of this paper were first obtained in October 2015 and first appeared online in December 2015. In March 2016, Nesterov and Stich have reproduced our results in a technical report [22].

# 1 Introduction

First-order methods have received extensive attention in the past two decades due to their ability to handle large-scale optimization problems. Recently, the development of *coordinate* versions of first-order methods have pushed their running times even faster. As a notable example, the state-of-the-art algorithm for empirical risk minimization (ERM) problems, up to a primal-dual transformation, is precisely accelerated coordinate descent [14].

In this paper, we consider the following unconstrained minimization problem[1]

$$\min_{x \in \mathbb{R}^n} f(x) \qquad (1.1)$$

where the objective $f \colon \mathbb{R}^n \to \mathbb{R}$ is continuously differentiable and convex. Below, we assume that $f(\cdot)$ is $L_i$-smooth with respect to its $i$-th coordinate.

Informally, coordinate smoothness means for each input $x$, if we add its $i$-th coordinate by at most $\delta$, the corresponding coordinate gradient $\nabla_i f(x + \delta \mathbf{e}_i)$ differs from $\nabla_i f(x)$ by at most $L_i$ times $|\delta|$. Under this definition, the larger $L_i$ is, the less smooth $f$ is along the $\mathbf{e}_i$ direction and therefore the harder it is to minimize $f$ along the $\mathbf{e}_i$ direction.[2] Intuitively, this implies we should spend more energy (i.e., assign more sampling probability) on coordinates with larger $L_i$. However, it was unclear what the best design is for such a distribution. In this paper, we present a clean and novel non-uniform sampling method which gives a faster convergence rate. Before going into the details, we first draw a distinction between non-accelerated and accelerated coordinate descent methods.

**Non-Accelerated vs. Accelerated Methods.** For smooth convex minimization, many first-order methods converge at a rate $1/\varepsilon$ to obtain an additive error $\varepsilon > 0$. In 1983, Nesterov demonstrated that a better and optimal rate $1/\sqrt{\varepsilon}$ can be obtained using his seminal accelerated gradient descent method. [18]

For this reason, people refer to methods converging at rate $1/\varepsilon$ as *non-accelerated* first-order methods, while those at rate $1/\sqrt{\varepsilon}$ as *accelerated* first-order methods. Similarly, when the objective $f(\cdot)$ is known to be strongly convex with parameter $\sigma > 0$, non-accelerated methods converge at a rate inversely proportional to $\sigma$, while accelerated ones converge at a rate inversely proportional to $\sqrt{\sigma}$. Although being much faster, accelerated first-order methods are also much more involved to design, see some recent attempts for designing accelerated methods in conceptually simpler manners [3, 6, 24, 35].

Such a distinction continues to hold on the *coordinate-gradient* setting. A coordinate descent method iteratively selects a coordinate $i \in [n]$ at random, and updates the iterate $x$ according to its coordinate gradient $\nabla_i f(x)$. As we shall see later, designing good sampling probabilities is well-studied for *non-accelerated* coordinate descent. In contrast, less is known in the more challenging accelerated regime, and we hope our work fills this gap.

We begin describing our result and compare it to the literature in the Euclidean norm case.

## 1.1 The Standard Euclidean Norm Case

In the non-accelerated world, in 2012, Nesterov [21] proposed a coordinate descent method called `RCDM`, which is a simple adaption of the full gradient descent method (see for instance the textbook [19]). At each iteration, `RCDM` selects a coordinate $i$ with probability proportional to $L_i$, and performs update $x' \leftarrow x - \frac{1}{L_i} \nabla_i f(x)$. The number of iterations required to reach an $\varepsilon$ error, denoted by $T$ in this paper, satisfies $T = O(\frac{\sum_i L_i}{\varepsilon} \|x_0 - x^*\|^2)$ for RCDM. Here, we denote by $x_0$ the starting vector, $x^*$ the minimizer of $f$, and $\|\cdot\|$ the $\ell_2$ Euclidean norm.

---

[1]The results of this paper generalize to the so-called proximal case that is to allow an additional separable term $\psi(x) \stackrel{\text{def}}{=} \sum_{i=1}^{n} \psi_i(x_i)$ to be added. The proofs require some non-trivial changes so we refrain from doing so in this version of the paper.

[2]For instance, if the $i$-th coordinate is selected, most coordinate-descent methods are only capable of performing an update $x' \leftarrow x - \frac{1}{L_i} \nabla_i f(x)$ with step length inversely proportional to $L_i$.

| **Paper** | Euclidean $\beta = 0$ Case | | $\beta \in [0,1]$ Case | | $\beta = 1$ Case | |
|---|---|---|---|---|---|---|
| | strongly convex | non-strongly convex | strongly convex | non-strongly convex | strongly convex | non-strongly convex |
| RCDM [21], | $\frac{\sum_i L_i}{\sigma} \log \frac{1}{\varepsilon}$ | $\frac{\sum_i L_i}{\varepsilon} \|x_0 - x^*\|^2$ | $\frac{S_{1-\beta}}{\sigma_\beta} \log \frac{1}{\varepsilon}$ | $\frac{S_{1-\beta}}{\varepsilon} \|x_0 - x^*\|^2_{\mathsf{L}_\beta}$ | $\frac{n}{\sigma_1} \log \frac{1}{\varepsilon}$ | $\frac{n}{\varepsilon} \|x_0 - x^*\|^2_{\mathsf{L}_1}$ |
| APCG [14], RBCD [15], Nesterov [21], APPROX [10] | - | - | - | - | $\frac{n}{\sqrt{\sigma_1}} \log \frac{1}{\varepsilon}$ | $\frac{n}{\sqrt{\varepsilon}} \|x_0 - x^*\|_{\mathsf{L}_1}$ |
| ACDM [13] | $\frac{\sqrt{n \sum_i L_i}}{\sqrt{\sigma}} \log \frac{1}{\varepsilon}$ | $\frac{\sqrt{n \sum_i L_i}}{\sqrt{\varepsilon}/\log \frac{1}{\varepsilon}} \|x_0 - x^*\|$ | $\frac{\sqrt{nS_{1-\beta}}}{\sqrt{\sigma_\beta}} \log \frac{1}{\varepsilon}$ | $\frac{\sqrt{nS_{1-\beta}}}{\sqrt{\varepsilon}/\log \frac{1}{\varepsilon}} \|x_0 - x^*\|_{\mathsf{L}_\beta}$ | $\frac{\sqrt{n}}{\sqrt{\sigma_1}} \log \frac{1}{\varepsilon}$ | $\frac{n}{\sqrt{\varepsilon}/\log \frac{1}{\varepsilon}} \|x_0 - x^*\|_{\mathsf{L}_1}$ |
| this paper | $\frac{\sum_i \sqrt{L_i}}{\sqrt{\sigma}} \log \frac{1}{\varepsilon}$ | $\frac{\sum_i \sqrt{L_i}}{\sqrt{\varepsilon}} \|x_0 - x^*\|$ | $\frac{S_{(1-\beta)/2}}{\sqrt{\sigma_\beta}} \log \frac{1}{\varepsilon}$ | $\frac{S_{(1-\beta)/2}}{\sqrt{\varepsilon}} \|x_0 - x^*\|_{\mathsf{L}_\beta}$ | $\frac{\sqrt{n}}{\sqrt{\sigma_1}} \log \frac{1}{\varepsilon}$ | $\frac{n}{\sqrt{\varepsilon}} \|x_0 - x^*\|_{\mathsf{L}_1}$ |

Table 1: Comparisons among coordinate descent methods, where $S_\alpha \stackrel{\text{def}}{=} \sum_i L_i^\alpha$.

This convergence rate is usually compared to that of full gradient descent: if $L$ is the global smoothness parameter of $f(\cdot)$, then full gradient descent converges in $T = O(\frac{L}{\varepsilon} \|x_0 - x^*\|^2)$ iterations. Since $L_i$ is never larger than $L$, and performing a coordinate descent step is usually $n$ times faster than a full gradient step, RCDM performs faster than gradient descent in most applications.

In the same paper [21], Nesterov also demonstrated the possibility of performing *accelerated coordinate gradient descent* via a simple adaption of its full-gradient variant [18–20]. This has been later analyzed in full by Lee and Sidford [13], and they named this method *accelerated coordinate descent method (*ACDM*)*. ACDM converges the following number of iterations:

$$T = \begin{cases} \widetilde{O}\Big(\frac{\sqrt{n \sum_i L_i}}{\sqrt{\varepsilon}} \|x_0 - x^*\|\Big), & \text{when } f \text{ is convex}^3 \\ O\Big(\frac{\sqrt{n \sum_i L_i}}{\sqrt{\sigma}} \log \frac{1}{\varepsilon}\Big), & \text{when } f \text{ is } \sigma\text{-strongly convex} \end{cases}$$

ACDM is built upon the estimation sequence technique of Nesterov [18, 19, 21], and similar to RCDM, ACDM also selects each coordinate $i$ (essentially) with a probability proportional to $L_i$.[4] Since the analysis of Lee and Sidford is tight, it has been thought that the iteration bound $T$ is not improvable.

In this paper, with a different non-uniform sampling method, we develop a new accelerated coordinate descent method NU_ACDM that converges in $T$ iterations, where

$$T = \begin{cases} O\Big(\frac{\sum_i \sqrt{L_i}}{\sqrt{\varepsilon}} \|x_0 - x^*\|\Big), & \text{when } f \text{ is convex} \\ O\Big(\frac{\sum_i \sqrt{L_i}}{\sqrt{\sigma}} \log \frac{1}{\varepsilon}\Big), & \text{when } f \text{ is } \sigma\text{-strongly convex} \end{cases}$$

---

[3]In fact, Lee and Sidford did not include a version of ACDM that works for non-strongly convex objectives. However, using regularization and simple reduction, one can turn an iterative solver for strongly convex and smooth minimization to that for convex and smooth minimization, by replacing $\sigma$ with $\varepsilon/\|x_0 - x^*\|^2$. Such a reduction incurs a factor $\log(1/\varepsilon)$ which we hide with the $\widetilde{O}$ notation in this paper.

[4]More precisely, they select each coordinate $i$ with a probability proportional to $\max\{L_i, \frac{1}{n} \sum_j L_j\}$. As a consequence, each coordinate $i$ is selected with probability at least $\Omega(1/n)$. Lee and Sidford emphasized that using this sampling method, rather than choosing each $i$ directly with probability $L_i/(\sum_j L_j)$, is essential for ACDM to obtain the accelerated convergence rate.

Note that NU_ACDM is always faster than ACDM because $\sum_i \sqrt{L_i} \leq \sqrt{n \sum_i L_i}$. In the case when $(L_1, \ldots, L_n)$ is non-uniform, our method runs faster by a factor up to $\sqrt{n}$.[5] In our sampling step, we select each coordinate $i$ with probability exactly proportional to $\sqrt{L_i}$, rather than (roughly) proportional to $L_i$. Thus, we need a different analysis from ACDM [13], and also avoid the more complicated estimation sequence analysis.

## 1.2 The General $\mathsf{L}_\beta$-Norm Case

Define the $\mathsf{L}_\beta$ norm $\|y\|_{\mathsf{L}_\beta}^2 \stackrel{\text{def}}{=} \sum_i L_i^\beta \cdot y_i^2$ for $\beta \in [0, 1]$. Many accelerated coordinate descent methods provide convergence guarantees with respect to the $\mathsf{L}_1$ norm [10, 15] or the $\mathsf{L}_\beta$ norm [13, 14, 21].

For instance, RCDM takes $\beta$ as an input, and converges in $T = O\big(\frac{S_{1-\beta}}{\varepsilon}\|x_0 - x^*\|_{\mathsf{L}_\beta}^2\big)$ iterations if one samples each coordinate $i$ with probability $L_i^{1-\beta}/S_{1-\beta}$, where $S_\alpha \stackrel{\text{def}}{=} \sum_i L_i^\alpha$. In [13], Lee and Sidford showed that their ACDM converges in $T$ iterations with the same sampling probabilities $L_i^{1-\beta}/S_{1-\beta}$, where

$$T = \begin{cases} \widetilde{O}\big(\frac{\sqrt{nS_{1-\beta}}}{\sqrt{\varepsilon}}\|x_0 - x^*\|_{\mathsf{L}_\beta}^2\big), & \text{when } f \text{ is convex} \\ O\big(\frac{\sqrt{nS_{1-\beta}}}{\sqrt{\sigma_\beta}} \log \frac{1}{\varepsilon}\big), & \text{when } f \text{ is } \sigma_\beta\text{-strongly convex w.r.t. the } \mathsf{L}_\beta \text{ norm} \end{cases}$$

This is always faster than RCDM. Note that, in the special case of $\beta = 1$ (and thus using uniform sampling probabilities), this same convergence result is also obtained by Nesterov [21], APCG [14], RBCD [15], and APPROX [10]. (See Table 1.) The results with respect to the $\mathsf{L}_1$ norm are not very interesting. One can scale each coordinate by a factor $1/L_i$ and apply an existing uniform-sampling coordinate-descent theorem to obtain the same result.

Our method NU_ACDM improves this convergence to

$$T = \begin{cases} O\big(\frac{S_{(1-\beta)/2}}{\sqrt{\varepsilon}}\|x_0 - x^*\|_{\mathsf{L}_\beta}\big), & \text{when } f \text{ is convex} \\ O\big(\frac{S_{(1-\beta)/2}}{\sqrt{\sigma_\beta}} \log \frac{1}{\varepsilon}\big), & \text{when } f \text{ is } \sigma_\beta\text{-strongly convex w.r.t. the } \mathsf{L}_\beta \text{ norm} \end{cases}$$

Since $\sqrt{S_{1-\beta}} \leq S_{(1-\beta)/2} \leq \sqrt{nS_{1-\beta}}$, our method is faster than ACDM by a factor up to $\sqrt{n}$. Our improvement is again due to the new choice of sampling probabilities —we select each coordinate $i$ with probability $L_i^{(1-\beta)/2}/S_{(1-\beta)/2}$ which is different from RCDM or ACDM— as well as our new proof that avoids the use of estimation sequence.

**Remark 1.1.** For the strongly convex case, convergence results with respect to Euclidean norms are usually more relevant to applications: for instance, the $\ell_2$ regularizer is the most common one used in machine learning applications and algorithms designed for the Euclidean norm should be used for a better performance.[6] However, in the non-strongly convex case, results with respect to different $\beta$ are in general incomparable. We include experiments in Section 7.3 to illustrate this.

---

[5]If $L_1 = \cdots = L_n$, we have $\sum_i \sqrt{L_i} = \sqrt{n \sum_i L_i}$. However, if $L_1 = 1$ while $L_2 = \cdots = L_n = 0$, we have $\sum_i \sqrt{L_i} = \frac{1}{\sqrt{n}} \cdot \sqrt{n \sum_i L_i}$.

[6]In contrast, consider an objective $f(x)$ equipped with a regularizer $\frac{\sigma}{2}\|x\|^2$. Such an objective is also strongly convex with respect to the $\mathsf{L}_\beta$ norm with parameter $\min_i L_i^{-\beta}$. If one applies an algorithm designed for the $\mathsf{L}_\beta$ norm using this parameter, the convergence would be much worse than the first column of Table 1.

## 2 Applications

**Empirical Risk Minimization.** A cornerstone problem in machine learning is empirical risk minimization (ERM). Let $a_1, \ldots, a_n \in \mathbb{R}^d$ be the feature vectors of $n$ data samples, $\phi_1, \ldots, \phi_n \colon \mathbb{R} \to \mathbb{R}$ be a sequence of convex loss functions, and $r \colon \mathbb{R}^d \to \mathbb{R}$ be a convex function (often known as a regularizer). The goal of ERM problem is to solve the following *primal* convex problem:

$$\min_{w \in \mathbb{R}^d} P(w) \stackrel{\text{def}}{=} \frac{1}{n} \sum_{i=1}^{n} \phi_i\big(\langle a_i, w \rangle\big) + r(w). \tag{2.1}$$

This includes a family of important problems such as SVM, Lasso, ridge regression, and logistic regression. Lin, Lu and Xiao [14] showed that the above minimization problem is equivalent to the following dual one:

$$\min_{y \in \mathbb{R}^n} D(y) \stackrel{\text{def}}{=} \frac{1}{n} \sum_{i=1}^{n} \phi_i^*(y_i) + r^*\left(-\frac{1}{n} \sum_{i=1}^{n} y_i a_i\right), \tag{2.2}$$

where $\phi_i^*$ and $r^*$ are respectively the Fenchel conjugate function of $\phi_i$ and $r$.[7] Most importantly, if properly preprocessed, $D(y)$ can be shown to be coordinate-wise smooth and therefore accelerated coordinate descent methods can be applied to minimize $D(y)$. This approach leads to algorithm APCG, which matches the best known worst-case running time on solving (2.1) up to a logarithmic factor.[8]

However, by taking a closer look, the coordinate smoothness parameters $L_1, \ldots, L_n$ of $D(y)$ are data dependent. Indeed, $L_i$ is roughly proportional to the Euclidean norm square of the $i$-th feature vector. Therefore, we can apply NU_ACDM in this paper to improve the running time obtained by APCG or AccSDCA. This is done in Section 7.

Note that each iteration of NU_ACDM selects a feature vector with a probability (roughly) proportional to its Euclidean norm. This is very different from the recent work of Zhao and Zhang [38], where they observed that for SDCA [32], a non-accelerated method, feature vectors should be sampled with probabilities proportional to their Euclidean norm *squares*. If one also uses the squared norms in the accelerated setting, he will only get a running time similar to ACDM, and therefore worse than our NU_ACDM.

We also mention one recent result that uses our NU_ACDM to develop faster ERM methods by exploiting the clustering structure of the dataset [4].

**Solving Linear Systems.** Consider a linear system $Ax = b$ for some full row rank matrix $A \in \mathbb{R}^{m \times n}$ where $m \geq n$. Denoting $a_i \in \mathbb{R}^n$ as the $i$-th row vector of matrix $A$, the celebrated Kaczmarz method [12] iteratively picks one of the row vectors $a_i$ and computes

$$x_{k+1} \leftarrow x_k + \frac{b_i - \langle a_i, x_k \rangle}{\|a_i\|^2} a_i \ .$$

Although many deterministic schemes have been proposed regarding how to select row vectors, many of them are difficult to analyze or compare. In a breakthrough paper, Strohmer and Vershynin [34] analyzed a randomized scheme and proved that:

---

[7]The conjugate of $r(x)$ is $r^*(y) \stackrel{\text{def}}{=} \max_w \{y^T w - r(w)\}$.

[8]Accelerated algorithms for solving (2.1) were first obtained by AccSDCA [33], and more recently improved by Katyusha [1].

**Theorem 2.1** (Randomized Kaczmarz [34]). *If one samples row $i$ with probability proportional to $\|a_i\|^2$ in each iteration, then the Kaczmarz method produces an $\varepsilon$-approximate solution of $Ax = b$[9] in $O(\|A^{-1}\|_2^2 \cdot \|A\|_F^2 \cdot \log \frac{1}{\varepsilon})$ iterations, and each iteration costs a running time $O(n)$.*

Above, $x^*$ is the solution to $Ax = b$, $A^{-1}$ is the left inverse, $\|A^{-1}\|_2$ is one divided by the smallest non-zero singular value of $A$, and $\|A\|_F = (\sum_{ij} a_{ij}^2)^{1/2}$ is the Frobenius norm.

Randomized Kaczmarz can be viewed as coordinate descent [11, 13, 17], and therefore `ACDM` applies here and gives a faster running time:

**Theorem 2.2** (`ACDM` on Kaczmarz [13]). *The `ACDM` method samples row $i$ with probability proportional to $\max\{\|a_i\|^2, \frac{\|A\|_F^2}{m}\}$ at each iteration, and produces an $\varepsilon$-approximate solution to $Ax = b$ in $O(\sqrt{m}\|A^{-1}\|_2 \cdot \|A\|_F \cdot \log \frac{1}{\varepsilon})$ iterations, and each iteration costs a running time $O(n)$.*

To obtain the above result, Lee and Sidford rewrote the problem of solving $Ax = b$ as an $m$-variate quadratic minimization problem

$$\min_{y \in \mathbb{R}^m} \left\{ f(y) \stackrel{\text{def}}{=} \frac{1}{2} \|A^T y\|^2 - \langle b, y \rangle \right\} .$$

The coordinate smoothness of $f$ is $L_i = \|a_i\|^2$ for every $i \in [m]$, and the strong convexity of $f$ can be deduced as $\sigma = \|A^{-1}\|_2^{-2}$.[10] For this reason, if we apply `NU_ACDM` instead of `ACDM`, we immediately get a faster algorithm:

**Theorem 2.3** (`NU_ACDM` on Kaczmarz). *The `NU_ACDM` method samples row $i$ with probability proportional to $\|a_i\|$ at each iteration, and produces an $\varepsilon$-approximate solution to $Ax = b$ in $O(\|A^{-1}\|_2 \cdot \|A\|_{2,1} \cdot \log \frac{1}{\varepsilon})$ iterations, and each iteration costs a running time $O(n)$.*

Above, $\|A\|_{2,1} \stackrel{\text{def}}{=} \sum_{j=1}^m \left( \sum_{i=1}^n |a_{ij}|^2 \right)^{1/2}$ is the matrix $L_{2,1}$ norm. Since it satisfies $\|A\|_F \leq \|A\|_{2,1} \leq \sqrt{m}\|A\|_F$, our method is always faster than `ACDM`, and can be faster by a factor up to $\sqrt{m}$ that depends on the problem structure. We provide empirical evaluation on this in Section 8.

We remark here that when $A$ is a positive semidefinite square matrix, the celebrated conjugate gradient method solves linear systems $Ax = b$ efficiently, especially when the eigenvalues of $A$ are clustered. The authors of `ACDM` have reported that their method can run faster than conjugate gradient in several interesting cases, and our `NU_ACDM` obviously provides further speedups to such cases as well.

## 3  Other Related Work

People have considered selecting coordinates non-uniformly from other perspectives. For example, Nutini et al. [23] compared the random coordinate selection rule with the Gauss-Southwell rule, and proved that except in the extreme cases, Gauss-Southwell rule is faster. Needell et al. [17] proposed a non-uniform sampling for stochastic gradient descent, and made a connection to the randomized Kaczmarz algorithm. Qu et al. [26] gave an algorithm which supports arbitrary sampling on dual variables. Csiba et al. [7] showed that one can adaptively choose a probability distribution over the dual variables that depends on the "dual residues". All of the works cited above are for non-accelerated settings, while this paper focuses on designing fast accelerated method. Note that Qu and Richtárik [25] provided a unified analysis for both accelerated and non-accelerated coordinate

---

[9]That is, a vector $x$ satisfying $\mathbb{E}[\|x_k - x^*\|^2] \leq \varepsilon \|x_0 - x^*\|^2$.
[10]One has to in fact consider the strong convexity of $f$ in the space orthogonal to the null space $\{y \in \mathbb{R}^m | A^T y = 0\}$. We recommend interested readers to see Section 5.2 of [13] for details.

**Algorithm 1** NU_ACDM($\beta, f, x_0, T$)

---
**Input:** $\beta \in [0, 1]$;
  $f$ a convex function that is coordinate-wise smooth with parameters $(L_1, \ldots, L_n)$, and $\sigma_\beta$-strongly convex with respect to $\|\cdot\|_{\mathsf{L}_\beta}$ for some $\beta \in [0, 1]$;
  $x_0$ some initial point; and
  $T$ the number of iterations.
**Output:** $y_T$ such that $\mathbb{E}[f(y_T)] - f(x^*) \leq O((1-\tau)^T) \cdot (f(x_0) - f(x^*))$.

1: $\alpha \leftarrow (1-\beta)/2$, $S_\alpha \leftarrow \sum_{i=1}^n L_i^\alpha$.
2: $p_i \leftarrow \frac{L_i^\alpha}{S_\alpha}$ for each $i \in [n]$. $\quad\diamond\ \sum_i p_i = 1$ so $\{p_i\}_i$ forms a distribution over $[n]$
3: $\tau \leftarrow \frac{2}{1+\sqrt{4S_\alpha^2/\sigma_\beta+1}}$, $\eta \leftarrow \frac{1}{\tau S_\alpha^2}$. $\quad\diamond\ \tau = O(\frac{\sqrt{\sigma_\beta}}{S_\alpha})$ and $\eta = O(\frac{1}{\sqrt{\sigma_\beta} S_\alpha})$
4: $y_0 \leftarrow x_0, \quad z_0 \leftarrow x_0$.
5: **for** $k \leftarrow 0$ **to** $T-1$ **do**
6: $\quad x_{k+1} \leftarrow \tau z_k + (1-\tau) y_k$.
7: $\quad$ Sample $i$ from $\{1, \cdots, n\}$ based on $p = (p_1, \cdots, p_n)$.
8: $\quad y_{k+1} \leftarrow y_{k+1}^{(i)} \stackrel{\mathrm{def}}{=} x_{k+1} - \frac{1}{L_i} \nabla_i f(x_{k+1})$
9: $\quad z_{k+1} \leftarrow z_{k+1}^{(i)} \stackrel{\mathrm{def}}{=} \frac{1}{1+\eta\sigma_\beta} \big( z_k + \eta\sigma_\beta x_{k+1} - \frac{\eta}{p_i L_i^\beta} \nabla_i f(x_{k+1}) \big)$
10: **end for**
11: **return** $y_T$.

---

descent methods with what they call "arbitrary sampling" in the non-strongly convex case. Our work can be seen as a continuation of that work, in that we instead focus on a particular class of sampling probabilities, for which we derive provably better convergence complexity bounds than prior results both for strongly-convex and non-strongly convex cases. In the non-strongly convex case, our results can be infered from the general results in [25].

This paper aims at improving the intrinsic convergence rate for coordinate descent, and thus chooses not to cover some of the other important features.

- One such feature is to support composite objective function of the form $\min_x \{f(x) + \sum_i \psi_i(x_i)\}$ for simple and proper convex proximal functions $\psi_i(\cdot)$. Many coordinate descent methods support this feature [9, 16, 28–31]. The recent results of RBCD [15] and APCG [14] made an important step towards supporting proximal functions also for accelerated coordinate descent. Our results apply to such setting as well, but since this will significantly complicate the notations and the proofs, we refrain from doing so in this paper.

- Another important feature is to support parallelization in coordinate descent. The interesting sequence of work [5, 27, 36] manages to prove that, very informally, if the function satisfies some nice property then the coordinate updates can be performed in parallel. This type of result has been recently introduced to accelerated coordinate descent as well [28].

For empirical risk minimization problems, we emphasize that accelerated coordinate descent (such as APCG [14]) is not the only tool to obtain the fastest running time. The first method with this running time is AccSDCA by Shalev-Shwartz and Zhang [33], and the best known running time is by Katyusha [1]. However, these methods do not solve (1.1) through coordinate gradients.

# 4 Notations

Let $x^*$ be an arbitrary minimizer of $f(x)$ and we are interested in finding a vector $x$ satisfying $f(x) - f(x^*) \leq \varepsilon$ for an accuracy parameter $\varepsilon > 0$. We use $\|\cdot\|$ to denote the Euclidean norm and $\mathbf{e}_i \in \mathbb{R}^n$ the $i$-th unit vector. We denote by $\nabla f(x)$ the full gradient of $f$ at point $x \in \mathbb{R}^n$, and by $\nabla_i f(x)$ the $i$-th coordinate gradient. With a slight abuse of notation, we view $\nabla_i f(x)$ both as a scaler in $\mathbb{R}$ and as a singleton vector in $\mathbb{R}^n$.

**Definition 4.1.** *We say that $f$ is $L$-smooth if $\forall x, y \in \mathbb{R}^n$, it satisfies $\|\nabla f(x) - \nabla f(y)\| \leq L\|x-y\|$.*
*We say that $f$ is $\sigma$-strongly convex (with respect to the Euclidean norm) if $\forall x, y \in \mathbb{R}^n$, it satisfies $f(y) \geq f(x) + \langle \nabla f(x), y - x \rangle + \frac{\sigma}{2}\|x - y\|^2$.*

**Definition 4.2.** *$f$ is coordinate-wise smooth with parameters $(L_1, L_2, \ldots, L_n)$, if for all $x \in \mathbb{R}^n$, $\delta > 0$, $i \in [n]$:*
$$|\nabla_i f(x + \delta \mathbf{e}_i) - \nabla_i f(x)| \leq L_i \cdot \delta .$$

We remark here that the same result of this paper continues to hold if one replaces coordinate gradients with the more general block-wise gradients.[11] We adopt this simpler notion for the ease of presenting this paper.

Following the notations of prior work [13, 21], we make the following definitions

**Definition 4.3.** *Given $\alpha, \beta \in [0, 1]$, define*
$$S_\alpha \stackrel{\text{def}}{=} \sum_{i=1}^n L_i^\alpha, \qquad \|x\|_{\mathsf{L}_\beta} \stackrel{\text{def}}{=} \sum_{i=1}^n x_i^2 \cdot L_i^\beta, \quad \text{and} \quad \langle x, y \rangle_{\mathsf{L}_\beta} \stackrel{\text{def}}{=} \sum_{i=1}^n x_i y_i \cdot L_i^\beta .$$

*Also, define $\sigma_\beta$ to be the strong convexity parameter of $f(\cdot)$ with respect to the $\|\cdot\|_{\mathsf{L}_\beta}$ norm. That is, it satisfies $f(y) \geq f(x) + \langle \nabla f(x), y - x \rangle + \frac{\sigma_\beta}{2}\|x - y\|_{\mathsf{L}_\beta}^2$ for all $x, y \in \mathbb{R}^n$.*

Clearly, if $f$ is $\sigma$ strongly convex then $\sigma_0 = \sigma$.

# 5 NUACDM in the Strongly Convex Case

We now propose our new method `NU_ACDM` to deal with strongly convex and smooth objectives. Suppose $f(\cdot)$ is coordinate-wise smooth with parameters $(L_1, \ldots, L_n)$ and $\sigma_\beta$-strongly convex with respect to $\|\cdot\|_{\mathsf{L}_\beta}$ for some $\beta \in [0, 1]$. At a first reading, one can simply consider $\beta = 0$ so $f$ is $\sigma_0$-strongly convex with respect to the traditional Euclidean norm. We choose to analyze the full parameter regime of $\beta$ to better compare us with known literatures.

As described in Algorithm 1, `NU_ACDM` begins with $x_0 = y_0 = z_0$ and iteratively computes the tuple $x_{k+1}, y_{k+1}, z_{k+1}$ from $x_k, y_k, z_k$. In iteration $k = 0, 1, \ldots, T-1$, we first compute $x_{k+1} \leftarrow \tau z_k + (1-\tau)y_k$ for some parameter $\tau \in [0, 1]$ (whose value will be specified later), and randomly select a coordinate $i \in [n]$ with probability $p_i = L_i^\alpha / S_\alpha$ where $\alpha \stackrel{\text{def}}{=} (1-\beta)/2$.

Whenever $i$ is selected at iteration $k$, we perform two updates $y_{k+1} \leftarrow x_{k+1} - \frac{1}{L_i} \nabla_i f(x_{k+1})$ and $z_{k+1} \leftarrow \frac{1}{1+\eta\sigma_\beta}\big(z_k + \eta\sigma_\beta x_{k+1} - \frac{\eta}{p_i L_i^\beta}\nabla_i f(x_{k+1})\big)$, both using the $i$-th coordinate gradient at point $x_{k+1}$. Here, $\eta > 0$ is the parameter that determines the step length of the second update; its choice will become clear in the analysis. Our main theorem in this section is as follows:

---
[11]That is, following for instance [14, 15], suppose that there is an $n \times n$ permutation matrix $U$ partitioned as $U = [U_1 \cdots U_p]$ where $U_i \in \mathbb{R}^{n \times n_i}$, and suppose that every vector $x \in \mathbb{R}^n$ is partitioned into $\{x_i \in \mathbb{R}^{n_i} \; i = 1, 2, \ldots, p\}$ satisfying that $x = \sum_{i=1}^p U_i x_i$. Next, one can define $\nabla_i f(x) \stackrel{\text{def}}{=} U_i^T \nabla f(x)$. Under these notations, we say that $f$ is block-wise smooth with parameters $(L_1, \ldots, L_p)$ if $\|\nabla_i f(x + U_i h_i) - \nabla_i f(x)\| \leq L_i \|h_i\|$ for every $i \in [p]$ and every $h_i \in \mathbb{R}^{n_i}$.

**Theorem 5.1.** *If $f(x)$ is coordinate-wise smooth with parameters $(L_1, \ldots, L_n)$, and $\sigma_\beta$-strongly convex with respect to $\|\cdot\|_{L_\beta}$ for some $\beta \in [0,1]$, then* `NU_ACDM`$(\beta, f, x_0, T)$ *produces an output $y_T$ satisfying*
$$\mathbb{E}[f(y_T)] - f(x^*) \leq O(1) \cdot (1-\tau)^T (f(x_0) - f(x^*)) \;,$$
*where* $\tau = \dfrac{2}{1+\sqrt{4S_{(1-\beta)/2}^2/\sigma_\beta + 1}} = \dfrac{1}{O\big(S_{(1-\beta)/2}/\sqrt{\sigma_\beta}\big)}$.

In particular, if $\beta = 0$ parameter $\tau$ becomes $\tau = \dfrac{1}{O\big(\sum_i \sqrt{L_i}/\sqrt{\sigma}\big)}$.

**Remark 5.2** (Per-Iteration Cost of `NU_ACDM`). Before analyzing the convergence, we emphasize that the computational cost of each iteration in `NU_ACDM` is dominated by (1) a coordinate gradient computation $\nabla_i f(\cdot)$ and (2) a constant number of $n$-dimensional vector computations such as $x_{k+1} \leftarrow \tau z_k + (1-\tau) y_k$. In many real-life applications, the complexity of (2) dominates that of (1). For this reason, many authors have discussed about how to carefully modify an accelerated algorithm to avoid $n$-dimensional vector computations. We refer interested readers to Section 4 of [14], and simply claim here without proof that our `NU_ACDM` admits the same modification.

## 5.1 Proof Outline

Our proof is different from the estimation sequence analysis used in `ACDM`, `RCDM`, or `APCG`, but follows from the linear-coupling framework of [3].

We use the superscript $^{(i)}$ on $y_{k+1}^{(i)}$ and $z_{k+1}^{(i)}$ to emphasize that the value depends on the choice of $i$. Therefore, $y_{k+1}$ equals $y_{k+1}^{(i)}$ with probability $p_i$, and similarly for $z$.

At each iteration $k$, the coordinate-wise smoothness classically gives a guarantee on the objective decrease:

**Lemma 5.3.** $f(y_{k+1}^{(i)}) \leq f(x_{k+1}) - \frac{1}{2L_i}\|\nabla_i f(x_{k+1})\|^2$.

On the other hand, our update on $z$ can be written in the following minimization form, known as mirror-descent form in optimization literatures:
$$z_{k+1}^{(i)} = \min_z \left\{ \frac{1}{2}\|z - z_k\|_{L_\beta}^2 + \frac{\eta}{p_i}\langle \nabla_i f(x_{k+1}), z\rangle + \frac{\eta \sigma_\beta}{2}\|z - x_{k+1}\|_{L_\beta}^2 \right\} \;. \tag{5.1}$$

We then apply a probabilistic analysis to deduce the following guarantee on mirror descent:

**Lemma 5.4.** *For every $u \in \mathbb{R}^n$,*
$$\frac{\eta}{p_i}\langle \nabla_i f(x_{k+1}), z_{k+1}^{(i)} - u\rangle - \frac{\eta \sigma_\beta}{2}\|x_{k+1} - u\|_{L_\beta}^2$$
$$\leq -\frac{1}{2}\|z_k - z_{k+1}^{(i)}\|_{L_\beta}^2 + \frac{1}{2}\|z_k - u\|_{L_\beta}^2 - \frac{1+\eta\sigma_\beta}{2}\|z_{k+1}^{(i)} - u\|_{L_\beta}^2 \;,$$

Suppose at the moment that $y_k, z_k$ and $x_{k+1}$ are fixed, and the only source of randomness comes from the random choice of $i$ when computing $z_{k+1}$ and $y_{k+1}$. Then, the following inequality can be deduced as is a nature linear combination of the two lemmas above:

**Lemma 5.5.** *For every $u \in \mathbb{R}^n$,*
$$\eta\langle \nabla f(x_{k+1}), z_k - u\rangle - \frac{\eta\sigma_\beta}{2}\|u - x_{k+1}\|_{L_\beta}^2$$
$$\leq \eta^2 S_\alpha^2 \big(f(x_{k+1}) - \mathbb{E}_i[f(y_{k+1})]\big) + \frac{1}{2}\|z_k - u\|_{L_\beta}^2 - \frac{1+\eta\sigma_\beta}{2}\mathbb{E}_i\big[\|z_{k+1} - u\|_{L_\beta}^2\big]$$

By taking into account our choice $x_{k+1} = \tau z_k + (1-\tau)y_k$ and the convexity of $f(\cdot)$, we can deduce that (again assuming $y_k$ and $z_k$ are fixed and the only randomness comes from the choice of $i$ to compute $y_{k+1}$ and $z_{k+1}$):

**Lemma 5.6.**

$$0 \leq \frac{(1-\tau)\eta}{\tau}(f(y_k) - f(x^*)) - \frac{\eta}{\tau}\mathbb{E}_i[f(y_{k+1}) - f(x^*)] + \frac{1}{2}\|z_k - x^*\|^2_{\mathsf{L}_\beta} - \frac{1+\eta\sigma_\beta}{2}\mathbb{E}_i\big[\|z_{k+1} - x^*\|^2_{\mathsf{L}_\beta}\big]$$

Finally, we choose $\tau = \frac{2}{1+\sqrt{4S_\alpha^2/\sigma_\beta + 1}} \leq \frac{\sqrt{\sigma_\beta}}{S_\alpha}$ so as to ensure that $1 + \eta\sigma_\beta = \frac{1}{1-\tau}$.[12] Under these parameter choices, Lemma 5.6 can be re-written as

$$\mathbb{E}_i\Big[\big(f(y_{k+1}) - f(x^*)\big) + \frac{\tau}{2\eta(1-\tau)}\|z_{k+1} - x^*\|^2_{\mathsf{L}_\beta}\Big] \leq (1-\tau)\Big(\big(f(y_k) - f(x^*)\big) + \frac{\tau}{2\eta(1-\tau)}\|z_k - x^*\|^2_{\mathsf{L}_\beta}\Big)$$

Telescoping it for all iterations $k$, we conclude that

$$\mathbb{E}[f(y_T)] - f(x^*) \leq (1-\tau)^T\Big(f(y_0) - f(x^*) + \frac{\tau}{2\eta}\|z_0 - x^*\|^2_{\mathsf{L}_\beta}\Big) \leq O(1) \cdot (1-\tau)^T(f(x_0) - f(x^*)) \ .$$

where the last inequality is because (i) $x_0 = y_0 = z_0$, (ii) $O(\tau/\eta) = O(\tau^2 S_\alpha^2) = O(\sigma_\beta)$ and (iii) the strong convexity of $f(\cdot)$ which implies $f(x_0) - f(x^*) \geq \frac{\sigma_\beta}{2}\|x_0 - x^*\|^2_{\mathsf{L}_\beta}$. This finishes the proof of our Theorem 5.1.

## 5.2 Proofs of Missing Lemmas

**Lemma 5.3.** $f(y_{k+1}^{(i)}) \leq f(x_{k+1}) - \frac{1}{2L_i}\|\nabla_i f(x_{k+1})\|^2$.

*Proof.* Abbreviating $x_{k+1}$ by $x$ and $y_{k+1}^{(i)}$ by $y$, we have that $x$ and $y$ only different at coordinate $i$ so we deduce that

$$f(y) - f(x) = \int_{\tau=0}^1 \langle \nabla f(x + \tau(y-x)), y - x\rangle d\tau$$

$$= \langle \nabla f(x), y - x\rangle + \int_{\tau=0}^1 \langle \nabla f(x + \tau(y-x)) - \nabla f(x), y - x\rangle d\tau$$

$$= \langle \nabla f(x), y - x\rangle + \int_{\tau=0}^1 \big(\nabla_i f(x + \tau(y_i - x_i)\mathbf{e}_i) - \nabla_i f(x)\big) \cdot (y_i - x_i) d\tau$$

$$\leq \langle \nabla f(x), y - x\rangle + \int_{\tau=0}^1 \tau L_i(y_i - x_i) \cdot (y_i - x_i) d\tau = \langle \nabla f(x), y - x\rangle + \frac{L_i}{2}\|y - x\|^2 \ .$$

Above, the only inequality is due to the coordinate-wise smoothness of $f(\cdot)$ (see Definition 4.2). □

**Lemma 5.4.** *For every $u \in \mathbb{R}^n$,*

$$\frac{\eta}{p_i}\langle \nabla_i f(x_{k+1}), z_{k+1}^{(i)} - u\rangle - \frac{\eta\sigma_\beta}{2}\|x_{k+1} - u\|^2_{\mathsf{L}_\beta}$$

$$\leq -\frac{1}{2}\|z_k - z_{k+1}^{(i)}\|^2_{\mathsf{L}_\beta} + \frac{1}{2}\|z_k - u\|^2_{\mathsf{L}_\beta} - \frac{1+\eta\sigma_\beta}{2}\|z_{k+1}^{(i)} - u\|^2_{\mathsf{L}_\beta} \ ,$$

---

[12]The reason of $\tau \in [0,1]$ is as follows. By the definition of coordinate-wise smoothness and strong convexity we have for every $\beta$, $L_i \geq \sigma_\beta L_i^\beta$. This means $L_i^{1-\beta} \geq \sigma_\beta$. Since $\beta = 1-2\alpha$, we have the following fact: $S_\alpha^2 \geq \sum_{i=1}^n L_i^{2\alpha} \geq n\sigma_{1-2\alpha} = n\sigma_\beta$.

*Proof.* The minimality condition in (5.1) tells us that for every $u \in \mathbb{R}^n$,

$$0 = \left\langle \frac{\partial}{\partial z}\left(\frac{1}{2}\|z - z_k\|_{\mathsf{L}_\beta}^2 + \frac{\eta}{p_i}\langle \nabla_i f(x_{k+1}), z\rangle + \frac{\eta \sigma_\beta}{2}\|z - x_{k+1}\|_{\mathsf{L}_\beta}^2\right)\bigg|_{z=z_{k+1}^{(i)}}, z_{k+1}^{(i)} - u \right\rangle$$
$$= \langle z_{k+1}^{(i)} - z_k, z_{k+1}^{(i)} - u\rangle_{\mathsf{L}_\beta} + \frac{\eta}{p_i}\langle \nabla_i f(x_{k+1}), z_{k+1}^{(i)} - u\rangle + \eta\sigma_\beta \langle z_{k+1}^{(i)} - x_{k+1}, z_{k+1}^{(i)} - u\rangle_{\mathsf{L}_\beta} \quad (5.2)$$

Next, the three-point equality of Euclidean norms tells us that

$$\langle z_{k+1}^{(i)} - z_k, z_{k+1}^{(i)} - u\rangle_{\mathsf{L}_\beta} = \frac{1}{2}\|z_k - z_{k+1}^{(i)}\|_{\mathsf{L}_\beta}^2 - \frac{1}{2}\|z_k - u\|_{\mathsf{L}_\beta}^2 + \frac{1}{2}\|z_{k+1}^{(i)} - u\|_{\mathsf{L}_\beta}^2, \quad (5.3)$$

as well as that

$$\langle z_{k+1}^{(i)} - x_{k+1}, z_{k+1}^{(i)} - u\rangle_{\mathsf{L}_\beta} = \frac{1}{2}\|x_{k+1} - z_{k+1}^{(i)}\|_{\mathsf{L}_\beta}^2 - \frac{1}{2}\|x_{k+1} - u\|_{\mathsf{L}_\beta}^2 + \frac{1}{2}\|z_{k+1}^{(i)} - u\|_{\mathsf{L}_\beta}^2. \quad (5.4)$$

Plugging (5.3) and (5.4) back to (5.2), we arrive at the equality

$$\frac{\eta}{p_i}\langle \nabla_i f(x_{k+1}), z_{k+1}^{(i)} - u\rangle + \eta\left(\frac{\sigma_\beta}{2}\|x_{k+1} - z_{k+1}^{(i)}\|_{\mathsf{L}_\beta}^2 - \frac{\sigma_\beta}{2}\|x_{k+1} - u\|_{\mathsf{L}_\beta}^2\right)$$
$$= -\frac{1}{2}\|z_k - z_{k+1}^{(i)}\|_{\mathsf{L}_\beta}^2 + \frac{1}{2}\|z_k - u\|_{\mathsf{L}_\beta}^2 - \frac{1 + \eta\sigma_\beta}{2}\|z_{k+1}^{(i)} - u\|_{\mathsf{L}_\beta}^2,$$

thus finishing the proof. $\square$

**Lemma 5.5.** *For every $u \in \mathbb{R}^n$,*

$$\eta\langle \nabla f(x_{k+1}), z_k - u\rangle - \frac{\eta\sigma_\beta}{2}\|u - x_{k+1}\|_{\mathsf{L}_\beta}^2$$
$$\leq \eta^2 S_\alpha^2\big(f(x_{k+1}) - \mathbb{E}_i[f(y_{k+1})]\big) + \frac{1}{2}\|z_k - u\|_{\mathsf{L}_\beta}^2 - \frac{1 + \eta\sigma_\beta}{2}\mathbb{E}_i\big[\|z_{k+1} - u\|_{\mathsf{L}_\beta}^2\big]$$

*Proof.* Combining Lemma 5.3 and Lemma 5.4 we deduce that for each $i \in [n]$,

$$\frac{\eta}{p_i}\langle \nabla_i f(x_{k+1}), z_k - u\rangle - \frac{\eta\sigma_\beta}{2}\|x_{k+1} - u\|_{\mathsf{L}_\beta}^2$$
$$\overset{①}{\leq} \frac{\eta}{p_i}\langle \nabla_i f(x_{k+1}), z_k - z_{k+1}^{(i)}\rangle - \frac{1}{2}\|z_k - z_{k+1}^{(i)}\|_{\mathsf{L}_\beta}^2 + \frac{1}{2}\|z_k - u\|_{\mathsf{L}_\beta}^2 - \frac{1 + \eta\sigma_\beta}{2}\|z_{k+1}^{(i)} - u\|_{\mathsf{L}_\beta}^2$$
$$\overset{②}{\leq} \frac{\eta^2}{2p_i^2 L_i^\beta}\|\nabla_i f(x_{k+1})\|^2 + \frac{1}{2}\|z_k - u\|_{\mathsf{L}_\beta}^2 - \frac{1 + \eta\sigma_\beta}{2}\|z_{k+1}^{(i)} - u\|_{\mathsf{L}_\beta}^2$$
$$\overset{③}{\leq} \frac{\eta^2 L_i}{p_i^2 L_i^\beta}\big(f(x_{k+1}) - f(y_{k+1}^{(i)})\big) + \frac{1}{2}\|z_k - u\|_{\mathsf{L}_\beta}^2 - \frac{1 + \eta\sigma_\beta}{2}\|z_{k+1}^{(i)} - u\|_{\mathsf{L}_\beta}^2$$
$$\overset{④}{=} \eta^2 S_\alpha^2\big(f(x_{k+1}) - f(y_{k+1}^{(i)})\big) + \frac{1}{2}\|z_k - u\|_{\mathsf{L}_\beta}^2 - \frac{1 + \eta\sigma_\beta}{2}\|z_{k+1}^{(i)} - u\|_{\mathsf{L}_\beta}^2.$$

Above, ① uses Lemma 5.4, ② uses the Cauchy-Schwarz inequality,[13] ③ uses Lemma 5.3, and ④ uses the choice of $p_i = L_i^\alpha/S_\alpha$ and $\beta = 1 - 2\alpha$. As a result, taking into account the randomness of $i$, we

---

[13] More specifically, $\frac{\eta}{p_i}\langle \nabla_i f(x_{k+1}), z_k - z_{k+1}^{(i)}\rangle = \langle \frac{\eta}{p_i L_i^\beta}\nabla_i f(x_{k+1}), z_k - z_{k+1}^{(i)}\rangle_{\mathsf{L}_\beta} \leq \frac{\eta^2}{2p_i^2 L_i^{2\beta}}\|\nabla_i f(x_{k+1})\|_{\mathsf{L}_\beta}^2 + \frac{1}{2}\|z_k - z_{k+1}^{(i)}\|_{\mathsf{L}_\beta}^2 = \frac{\eta^2}{2p_i^2 L_i^\beta}\|\nabla_i f(x_{k+1})\|^2 + \frac{1}{2}\|z_k - z_{k+1}^{(i)}\|_{\mathsf{L}_\beta}^2$.

have

$$\eta \langle \nabla f(x_{k+1}), z_k - u \rangle - \frac{\eta \sigma_\beta}{2} \|u - x_{k+1}\|_{\mathsf{L}_\beta}^2$$

$$= \mathbb{E}_i \Big[ \frac{\eta}{p_i} \langle \nabla_i f(x_{k+1}), z_k - u \rangle - \frac{\eta \sigma_\beta}{2} \|u - x_{k+1}\|_{\mathsf{L}_\beta}^2 \Big]$$

$$\leq \mathbb{E}_i \Big[ \eta^2 S_\alpha^2 \big(f(x_{k+1}) - f(y_{k+1}^{(i)})\big) + \frac{1}{2} \|z_k - u\|_{\mathsf{L}_\beta}^2 - \frac{1 + \eta \sigma_\beta}{2} \|z_{k+1}^{(i)} - u\|_{\mathsf{L}_\beta}^2 \Big]$$

$$= \eta^2 S_\alpha^2 \big(f(x_{k+1}) - \mathbb{E}_i[f(y_{k+1})]\big) + \frac{1}{2} \|z_k - u\|_{\mathsf{L}_\beta}^2 - \frac{1 + \eta \sigma_\beta}{2} \mathbb{E}_i\big[\|z_{k+1} - u\|_{\mathsf{L}_\beta}^2\big] \ . \qquad \square$$

**Lemma 5.6.**

$$0 \leq \frac{(1-\tau)\eta}{\tau}(f(y_k) - f(x^*)) - \frac{\eta}{\tau}\mathbb{E}_i[f(y_{k+1}) - f(x^*)] + \frac{1}{2}\|z_k - x^*\|_{\mathsf{L}_\beta}^2 - \frac{1 + \eta \sigma_\beta}{2}\mathbb{E}_i\big[\|z_{k+1} - x^*\|_{\mathsf{L}_\beta}^2\big]$$

*Proof.*

$$\eta(f(x_{k+1}) - f(x^*))$$

$$\overset{\text{①}}{\leq} \eta \langle \nabla f(x_{k+1}), x_{k+1} - x^* \rangle - \frac{\eta \sigma_\beta}{2} \|x^* - x_{k+1}\|_{\mathsf{L}_\beta}^2$$

$$= \eta \langle \nabla f(x_{k+1}), x_{k+1} - z_k \rangle + \eta \langle \nabla f(x_{k+1}), z_k - x^* \rangle - \frac{\eta \sigma_\beta}{2} \|x^* - x_{k+1}\|_{\mathsf{L}_\beta}^2$$

$$\overset{\text{②}}{=} \frac{(1-\tau)\eta}{\tau} \langle \nabla f(x_{k+1}), y_k - x_{k+1} \rangle + \eta \langle \nabla f(x_{k+1}), z_k - x^* \rangle - \frac{\eta \sigma_\beta}{2} \|x^* - x_{k+1}\|_{\mathsf{L}_\beta}^2$$

$$\overset{\text{③}}{\leq} \frac{(1-\tau)\eta}{\tau}(f(y_k) - f(x_{k+1})) + \mathbb{E}_i\big[\eta^2 S_\alpha^2\big(f(x_{k+1}) - f(y_{k+1})\big) + \frac{1}{2}\|z_k - x^*\|_{\mathsf{L}_\beta}^2 - \frac{1 + \eta \sigma_\beta}{2}\|z_{k+1} - x^*\|_{\mathsf{L}_\beta}^2\big]$$

Above, ① is owing to the strong convexity of $f(\cdot)$ (see Definition 4.3), ② uses the fact that $x_{k+1} = \tau z_k + (1-\tau)y_k$, and ③ uses the convexity of $f(\cdot)$ as well as Lemma 5.5 with the choice of $u = x^*$. Recall $\eta = \frac{1}{\tau S_\alpha^2}$, we arrive at the desired inequality. $\qquad \square$

# 6 NUACDM in the Non-Strongly Convex Case

We now propose our non-uniform accelerated coordinate gradient method NU_ACDM$^{\mathrm{ns}}$ that deals with non-strongly convex and smooth objectives $f(\cdot)$. More precisely, consider the case when $f(\cdot)$ is coordinate-wise smooth with parameters $(L_1, \ldots, L_n)$.

As described in Algorithm 2, NU_ACDM$^{\mathrm{ns}}$ begins with $x_0 = y_0 = z_0$, and is parameterized by $\beta \in [0, 1]$. Following [13, 21], we shoot for obtaining a convergence result where the number of iterations $T$ is proportional to $\|x_0 - x^*\|_{\mathsf{L}_\beta}$.

Our NU_ACDM$^{\mathrm{ns}}$ iteratively compute the tuple $x_{k+1}, y_{k+1}, z_{k+1}$ from $x_k, y_k, z_k$. In iteration $k = 0, 1, \ldots, T-1$, we first compute $x_{k+1} \leftarrow \tau_k z_k + (1-\tau_k)y_k$ for some parameter $\tau_k \in [0, 1]$ (whose value will be specified and used at the end of this section), and randomly select a coordinate $i \in [n]$ with probability $p_i = L_i^\alpha / S_\alpha$ where $\alpha \overset{\mathrm{def}}{=} (1-\beta)/2$.

Whenever $i$ is selected at iteration $k$, we perform two updates $y_{k+1} \leftarrow x_{k+1} - \frac{1}{L_i} \nabla f_i(x_{k+1})$ and $z_{k+1} \leftarrow z_k - \frac{\eta_{k+1}}{p_i L_i^\beta} \nabla_i f(x_{k+1})$, both using the $i$-th coordinate gradient at point $x_{k+1}$. Here, $\eta_{k+1} > 0$ is the parameter that determines the step length of the second update; its choice will become clear at the end of this section.

We are now ready to state our main theorem in this section, and leave its proof to Section 6.1.

**Algorithm 2** NU_ACDM$^{\text{ns}}(\beta, f, x_0, T)$

**Input:** $\beta \in [0, 1]$;
  $f$ a convex function that is coordinate-wise smooth with parameters $(L_1, \ldots, L_n)$;
  $x_0$ some initial point; and
  $T$ the number of iterations.

**Output:** $y_T$ such that $\mathbb{E}[f(y_T)] - f(x^*) \leq 2\|x_0 - x^*\|^2_{\mathsf{L}_\beta} \cdot S^2_{(1-\beta)/2}/(T+1)^2$.

1: $\alpha \leftarrow (1-\beta)/2$, $S_\alpha \leftarrow \sum_{i=1}^n L_i^\alpha$.
2: $p_i \leftarrow \frac{L_i^\alpha}{S_\alpha}$ for each $i \in [n]$. ⋄ $\sum_i p_i = 1$ so $\{p_i\}_i$ forms a distribution over $[n]$
3: $y_0 \leftarrow x_0$, $z_0 \leftarrow x_0$.
4: **for** $k \leftarrow 0$ **to** $T - 1$ **do**
5: $\quad \eta_{k+1} \leftarrow \frac{k+2}{2S_\alpha^2}$, and $\tau_k \leftarrow \frac{1}{\eta_{k+1}S_\alpha^2} = \frac{2}{k+2}$.
6: $\quad x_{k+1} \leftarrow \tau_k z_k + (1 - \tau_k) y_k$.
7: $\quad$ Sample $i$ from $\{1, \cdots, n\}$ based on $p = (p_1, \cdots, p_n)$.
8: $\quad y_{k+1} \leftarrow y_{k+1}^{(i)} \stackrel{\text{def}}{=} x_{k+1} - \frac{1}{L_i}\nabla_i f(x_{k+1})$
9: $\quad z_{k+1} \leftarrow z_{k+1}^{(i)} \stackrel{\text{def}}{=} z_k - \frac{\eta_{k+1}}{p_i L_i^\beta}\nabla_i f(x_{k+1})$
10: **end for return** $y_T$.

---

**Theorem 6.1.** *If $f(x)$ is coordinate-wise smooth with parameters $(L_1, \ldots, L_n)$, and $\beta \in [0, 1]$ is a given parameter, the algorithm* NU_ACDM$^{\text{ns}}(\beta, f, x_0, T)$ *in Algorithm 2 produces an output $y_T$ satisfying*

$$\mathbb{E}[f(y_T)] - f(x^*) \leq \frac{2\|x_0 - x^*\|^2_{\mathsf{L}_\beta} \cdot S^2_{(1-\beta)/2}}{(T+1)^2} .$$

If $\beta = 0$, the above convergence gets simplified to

$$\mathbb{E}[f(y_T)] - f(x^*) \leq \frac{2\|x_0 - x^*\|^2 \cdot \left(\sum_i \sqrt{L_i}\right)^2}{(T+1)^2} .$$

### 6.1 Convergence Analysis

In this section we use the superscript $^{(i)}$ on $y_{k+1}^{(i)}$ and $z_{k+1}^{(i)}$ to emphasize that the value depends on the choice of $i$. Therefore, $y_{k+1}$ equals $y_{k+1}^{(i)}$ with probability $p_i$, and similarly for $z$.

At each iteration $k$, from Lemma 5.3 we know that the coordinate-wise smoothness directly yields the following guarantee on the objective decrease: $f(y_{k+1}^{(i)}) \leq f(x_{k+1}) - \frac{1}{2L_i}\|\nabla_i f(x_{k+1})\|^2$.

Next, since our update step on $z$ can be re-written in the minimization form

$$z_{k+1}^{(i)} = \min_z \left\{\frac{1}{2}\|z - z_k\|^2_{\mathsf{L}_\beta} + \frac{\eta_{k+1}}{p_i}\langle\nabla_i f(x_{k+1}), z\rangle\right\} , \quad (6.1)$$

we can apply the standard mirror-descent analysis and deduce that

**Lemma 6.2.** *For every $u \in \mathbb{R}^n$,*

$$\frac{\eta_{k+1}}{p_i}\langle\nabla_i f(x_{k+1}), z_{k+1}^{(i)} - u\rangle = -\frac{1}{2}\|z_k - z_{k+1}^{(i)}\|^2_{\mathsf{L}_\beta} + \frac{1}{2}\|z_k - u\|^2_{\mathsf{L}_\beta} - \frac{1}{2}\|z_{k+1}^{(i)} - u\|^2_{\mathsf{L}_\beta} ,$$

*Proof.* The minimality condition in (6.1) tells us that for every $u \in \mathbb{R}^n$,

$$0 = \left\langle \frac{\partial}{\partial z}\left(\frac{1}{2}\|z - z_k\|_{\mathsf{L}_\beta}^2 + \frac{\eta_{k+1}}{p_i}\langle \nabla_i f(x_{k+1}), z\rangle\right)\Big|_{z=z_{k+1}^{(i)}} , z_{k+1}^{(i)} - u \right\rangle$$

$$= \langle z_{k+1}^{(i)} - z_k, z_{k+1}^{(i)} - u\rangle_{\mathsf{L}_\beta} + \frac{\eta_{k+1}}{p_i}\langle \nabla_i f(x_{k+1}), z_{k+1}^{(i)} - u\rangle \quad (6.2)$$

Next, the three-point equality of Euclidean norms tells us that

$$\langle z_{k+1}^{(i)} - z_k, z_{k+1}^{(i)} - u\rangle_{\mathsf{L}_\beta} = \frac{1}{2}\|z_k - z_{k+1}^{(i)}\|_{\mathsf{L}_\beta}^2 - \frac{1}{2}\|z_k - u\|_{\mathsf{L}_\beta}^2 + \frac{1}{2}\|z_{k+1}^{(i)} - u\|_{\mathsf{L}_\beta}^2 . \quad (6.3)$$

Combining (6.3) and (6.2), we arrive at the desired equality. □

Suppose at the moment that $y_k, z_k$ and $x_{k+1}$ are fixed, and the only source of randomness comes from the random choice of $i$ when computing $z_{k+1}$ and $y_{k+1}$. We claim that the following inequality holds:

**Lemma 6.3.** *For every $u \in \mathbb{R}^n$,*

$$\eta_{k+1}\langle \nabla f(x_{k+1}), z_k - u\rangle \leq \eta_{k+1}^2 S_\alpha^2\big(f(x_{k+1}) - \mathbb{E}_i[f(y_{k+1})]\big) + \frac{1}{2}\|z_k - u\|_{\mathsf{L}_\beta}^2 - \frac{1}{2}\mathbb{E}_i\big[\|z_{k+1} - u\|_{\mathsf{L}_\beta}^2\big]$$

*Proof.* Combining Lemma 5.3 and Lemma 6.2 we deduce that for each $i \in [n]$,

$$\frac{\eta_{k+1}}{p_i}\langle \nabla_i f(x_{k+1}), z_k - u\rangle$$

$$\overset{\text{①}}{\leq} \frac{\eta_{k+1}}{p_i}\langle \nabla_i f(x_{k+1}), z_k - z_{k+1}^{(i)}\rangle - \frac{1}{2}\|z_k - z_{k+1}^{(i)}\|_{\mathsf{L}_\beta}^2 + \frac{1}{2}\|z_k - u\|_{\mathsf{L}_\beta}^2 - \frac{1}{2}\|z_{k+1}^{(i)} - u\|_{\mathsf{L}_\beta}^2$$

$$\overset{\text{②}}{\leq} \frac{\eta_{k+1}^2}{2p_i^2 L_i^\beta}\|\nabla_i f(x_{k+1})\|^2 + \frac{1}{2}\|z_k - u\|_{\mathsf{L}_\beta}^2 - \frac{1}{2}\|z_{k+1}^{(i)} - u\|_{\mathsf{L}_\beta}^2$$

$$\overset{\text{③}}{\leq} \frac{\eta_{k+1}^2 L_i}{p_i^2 L_i^\beta}\big(f(x_{k+1}) - f(y_{k+1}^{(i)})\big) + \frac{1}{2}\|z_k - u\|_{\mathsf{L}_\beta}^2 - \frac{1}{2}\|z_{k+1}^{(i)} - u\|_{\mathsf{L}_\beta}^2$$

$$\overset{\text{④}}{=} \eta_{k+1}^2 S_\alpha^2\big(f(x_{k+1}) - f(y_{k+1}^{(i)})\big) + \frac{1}{2}\|z_k - u\|_{\mathsf{L}_\beta}^2 - \frac{1}{2}\|z_{k+1}^{(i)} - u\|_{\mathsf{L}_\beta}^2 .$$

Above, ① uses Lemma 6.2, ② uses the Cauchy-Schwarz inequality,[14] ③ uses Lemma 5.3, and ④ uses the choice of $p_i = L_i^\alpha/S_\alpha$ and $\beta = 1 - 2\alpha$. As a result, taking into account the randomness of $i$, we have

$$\eta_{k+1}\langle \nabla f(x_{k+1}), z_k - u\rangle = \mathbb{E}_i\Big[\frac{\eta_{k+1}}{p_i}\langle \nabla_i f(x_{k+1}), z_k - u\rangle\Big]$$

$$\leq \mathbb{E}_i\big[\eta_{k+1}^2 S_\alpha^2\big(f(x_{k+1}) - f(y_{k+1}^{(i)})\big) + \frac{1}{2}\|z_k - u\|_{\mathsf{L}_\beta}^2 - \frac{1}{2}\|z_{k+1}^{(i)} - u\|_{\mathsf{L}_\beta}^2\big]$$

$$= \eta_{k+1}^2 S_\alpha^2\big(f(x_{k+1}) - \mathbb{E}_i[f(y_{k+1})]\big) + \frac{1}{2}\|z_k - u\|_{\mathsf{L}_\beta}^2 - \frac{1}{2}\mathbb{E}_i\big[\|z_{k+1} - u\|_{\mathsf{L}_\beta}^2\big] . \quad \square$$

At this moment, we take into account our choice of $x_{k+1} = \tau_k z_k + (1 - \tau_k)y_k$ and deduce that (again assuming that $y_k$ and $z_k$ are fixed while the only randomness comes from the choice of $i$ when computing $y_{k+1}$ and $z_{k+1}$):

---

[14]More specifically, $\frac{\eta_{k+1}}{p_i}\langle \nabla_i f(x_{k+1}), z_k - z_{k+1}^{(i)}\rangle = \langle \frac{\eta_{k+1}}{p_i L_i^\beta}\nabla_i f(x_{k+1}), z_k - z_{k+1}^{(i)}\rangle_{\mathsf{L}_\beta} \leq \frac{\eta_{k+1}^2}{2p_i^2 L_i^{2\beta}}\|\nabla_i f(x_{k+1})\|_{\mathsf{L}_\beta}^2 + \frac{1}{2}\|z_k - z_{k+1}^{(i)}\|_{\mathsf{L}_\beta}^2 = \frac{\eta_{k+1}^2}{2p_i^2 L_i^\beta}\|\nabla_i f(x_{k+1})\|^2 + \frac{1}{2}\|z_k - z_{k+1}^{(i)}\|_{\mathsf{L}_\beta}^2$.

**Lemma 6.4.** *For every $u \in \mathbb{R}^d$,*

$$0 \leq (\eta_{k+1}^2 S_\alpha^2 - \eta_{k+1})(f(y_k) - f(u)) - \eta_{k+1}^2 S_\alpha^2 \cdot \mathbb{E}_i\big[f(y_{k+1}) - f(u)\big] + \frac{1}{2}\|z_k - u\|_{\mathsf{L}_\beta}^2 - \frac{1}{2}\mathbb{E}_i\big[\|z_{k+1} - u\|_{\mathsf{L}_\beta}^2\big]$$

*Proof.*

$$\eta_{k+1}(f(x_{k+1}) - f(u)) \overset{\text{\textcircled{1}}}{\leq} \eta_{k+1}\langle \nabla f(x_{k+1}), x_{k+1} - u\rangle$$
$$= \eta_{k+1}\langle \nabla f(x_{k+1}), x_{k+1} - z_k\rangle + \eta_{k+1}\langle \nabla f(x_{k+1}), z_k - u\rangle$$
$$\overset{\text{\textcircled{2}}}{=} \frac{(1-\tau_k)\eta_{k+1}}{\tau_k}\langle \nabla f(x_{k+1}), y_k - x_{k+1}\rangle + \eta_{k+1}\langle \nabla f(x_{k+1}), z_k - u\rangle$$
$$\overset{\text{\textcircled{3}}}{\leq} \frac{(1-\tau_k)\eta_{k+1}}{\tau_k}(f(y_k) - f(x_{k+1})) + \mathbb{E}_i\big[\eta_{k+1}^2 S_\alpha^2\big(f(x_{k+1}) - f(y_{k+1})\big) + \frac{1}{2}\|z_k - u\|_{\mathsf{L}_\beta}^2 - \frac{1}{2}\|z_{k+1} - u\|_{\mathsf{L}_\beta}^2\big]$$

Above, ① is owing to the convexity of $f(\cdot)$, ② uses the fact that $x_{k+1} = \tau_k z_k + (1-\tau_k)y_k$, and ③ uses the convexity of $f(\cdot)$ as well as Lemma 6.3. As a consequence, by choosing $\tau_k = \frac{1}{\eta_{k+1}S_\alpha^2}$, we arrive at the desired inequality. □

Finally, we only need to set the sequence of $\eta_k$ so that $\eta_k^2 S_\alpha^2 \approx \eta_{k+1}^2 S_\alpha^2 - \eta_{k+1}$ as well as $\tau_k = 1/\eta_{k+1}S_\alpha^2 \in [0,1]$. For instance, we can let $\eta_k = \frac{k+1}{2S_\alpha^2}$ so that $\eta_k^2 S_\alpha^2 = \eta_{k+1}^2 S_\alpha^2 - \eta_{k+1} + \frac{1}{4S_\alpha^2}$. After telescoping Lemma 6.4 with $k = 0, 1, \ldots, T-1$ and setting $u = x^*$, we obtain that

$$\eta_T^2 S_\alpha^2 \mathbb{E}\big[f(y_T) - f(x^*)\big] + \sum_{k=1}^{T-1} \frac{1}{4S_\alpha^2}\mathbb{E}\big[f(y_k) - f(x^*)\big] \leq \frac{1}{2}\|z_0 - x^*\|_{\mathsf{L}_\beta}^2 - \frac{1}{2}\mathbb{E}\big[\|z_T - x^*\|_{\mathsf{L}_\beta}^2\big] \ .$$

Finally, since $x^*$ is the minimizer and satisfies $f(y_k) \geq f(x^*)$ and $z_0 = x_0$, we obtain

$$\frac{(T+1)^2}{4S_\alpha^4}S_\alpha^2\big(\mathbb{E}[f(y_T)] - f(x^*)\big) \leq \frac{1}{2}\|x_0 - x^*\|_{\mathsf{L}_\beta}^2 \ .$$

This finishes the proof of our Theorem 6.1.

# 7 Experiments on Empirical Risk Minimization

We perform experiments on ERM problems to confirm our theoretical improvements. We consider three datasets in this section: (1) class 1 of the `news20` dataset ($15,935$ samples and $62,061$ features), (2) the `w8a` dataset ($49,749$ samples and $300$ features), and (3) the `covtype` dataset ($581,012$ samples and $54$ features). All of them can be found on the LibSVM website [8], and contain examples that have non-uniform Euclidean norms (see Figure 1 for the distribution).

## 7.1 Experiments on Strongly Convex Objectives

Consider a regularized least-square problem which is problem (2.1) with $\phi_i(t) \overset{\text{def}}{=} \frac{1}{2}(t - l_i)^2$, where $l_i$ is the label for feature vector $a_i$. In the case when $r(w) = \frac{\lambda}{2}\|w\|_2^2$, this problem becomes *ridge regression*, and in the case when $r(w) = \lambda\|w\|_1$, it is known as *Lasso regression*.

Following (2.2), the equivalent dual formulation of regularized least square can be written as

$$\text{Dual:} \quad \min_{y \in \mathbb{R}^n}\left\{D(y) \overset{\text{def}}{=} \frac{1}{n}\sum_{i=1}^n \big(\frac{1}{2}y_i^2 + y_i \cdot l_i\big) + r^*\Big(-\frac{1}{n}\sum_{i=1}^n y_i a_i\Big)\right\} \ . \tag{7.1}$$

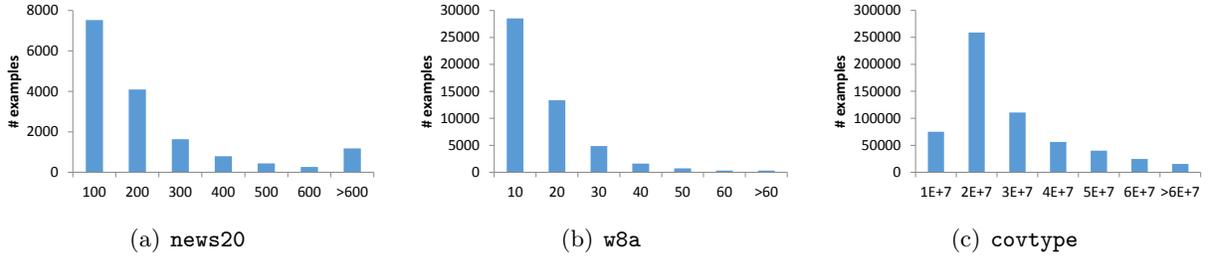

Figure 1: Distribution of $\|a_i\|_2^2$, the feature vector norm squares.

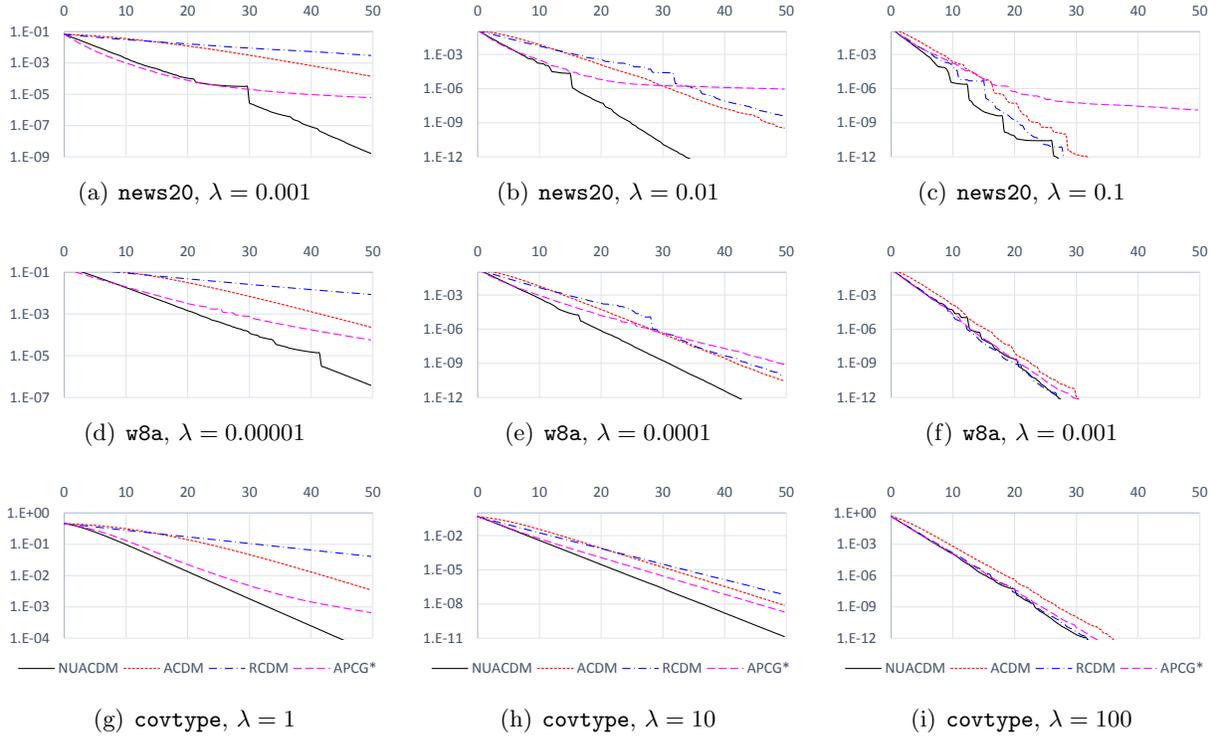

Figure 2: Performance Comparison for Ridge Regression. The $y$ axis represents the dual objective distance to minimum, and the $x$ axis represents the number of passes to the dataset.

Furthermore, $D(y)$ is $1/n$-strongly convex.

**Ridge Regression.** In ridge regression, we have $r(w) = \frac{\lambda}{2}\|w\|_2^2$ and accordingly $r^*(z) = \frac{1}{2\lambda}\|z\|_2^2$ in (7.1). It is not hard to verify that $D(y)$ is $L_i \stackrel{\text{def}}{=} \frac{1}{n} + \frac{1}{\lambda n^2}\|a_i\|_2^2$ smooth with respect to its $i$-th coordinate (and thus with respect to the $i$-th example). Therefore, the coordinate smoothness parameters are *non-uniform* if examples $a_1, \ldots, a_n$'s do not have the same Euclidean norms.

We can directly apply RCDM, ACDM and our NU_ACDM with $\beta = 0$ and $\sigma = 1/n$ to minimize (7.1). In principle, one can also apply APCG to minimize $D(y)$. However, since APCG is only designed for $\beta = 1$ and needs an unknown parameter $\sigma_1 > 0$ as input, we have tuned it for the fastest convergence in our experiments; whenever we do so, we denote it as APCG* in the diagrams.[15]

---
[15] We have chosen 14 values of $\sigma_1$ in a reasonable range, where the largest choice of $\sigma_1$ is $50,000$ times larger than the smallest choice. Our automated program will then make the final choice of $\sigma_1$ based on the convergence speed.

| | | | | | | |
|---|---|---|---|---|---|---|
| `news20` | $\lambda = 0.001$ | 1.56772 | $\lambda = 0.01$ | 1.30740 | $\lambda = 0.1$ | 1.05110 |
| `w8a` | $\lambda = 0.00001$ | 1.11060 | $\lambda = 0.0001$ | 1.04897 | $\lambda = 0.001$ | 1.00373 |
| `covtype` | $\lambda = 1$ | 1.04266 | $\lambda = 10$ | 1.02787 | $\lambda = 100$ | 1.00362 |

Table 2: The theoretical speed-up factor $\sqrt{n \sum_i L_i}/\left(\sum_i \sqrt{L_i}\right)$ of `NU_ACDM` over `ACDM` for the three datasets.

Our experimental results for ridge regression are in Figure 2. Note that theory predicts that `NU_ACDM` enjoys a speed-up factor of $\frac{\sqrt{n \sum_i L_i}}{\sum_i \sqrt{L_i}} \geq 1$ over `ACDM`, and we show this factor in Table 2. We make the following observations:

- Since $L_i = \frac{1}{n} + \frac{1}{\lambda n^2}\|a_i\|_2^2$, the smaller the regularization parameter $\lambda$ is, the more non-uniform the parameters $L_1, \ldots, L_n$ are. This is why the numbers in Table 2 are in decreasing order in each row. Our experiment confirms on this because we obtain the greatest improvements for the left 3 charts in Figure 2.

- `news20` has the most non-uniformity on the examples' Euclidean norms among the three datasets. Therefore, the first row Table 2 have the largest speed-up factors. Our experiment confirms on this because we obtain the greatest improvements in the top 3 charts in Figure 2.

- `APCG` performs quite poorly on dataset `news20` because it relies on the $\mathsf{L}_\beta$ norm strong convexity for $\beta = 1$, which is very different from the Euclidean norm strong convexity when the parameters $L_i$ are very non-uniform. We discuss the choice of $\beta$ in Section 7.3, and would like to point out that `APCG` performs very well for non strongly convex objectives, see Section 7.2.

Due to strong duality, our convergence speed-up on the dual objective also translates to that on the primal objective. See Figure 7 in the appendix for details.

**Lasso.** In the Lasso problem, we have $r(w) = \lambda\|w\|_1$ in the primal objective so the corresponding dual $D(y)$ in (7.1) is not smooth. Fortunately, standard regularization techniques suggest that in order to minimize the Lasso objective $P(w)$, it suffices to look at an alternative regularizer $r'(w) \stackrel{\text{def}}{=} \lambda\|w\|_1 + \frac{\lambda_2}{2}\|w\|_2^2$ and its corresponding objective $P'(w) \stackrel{\text{def}}{=} P(w) + \frac{\lambda_2}{2}\|w\|_2^2$.[16] Since for every $w$ it satisfies $|P(w) - P'(w)| \leq O(\lambda_2)$, one can specify a small enough parameter $\lambda_2 > 0$ and minimize $P'(w)$ instead. This auxiliary 2-norm regularizer introduces error to the objective, but allows the function $P'(w)$ to be minimized fast. Indeed, all known accelerated gradient methods to solve Lasso (such as `AccSDCA` [33], `APCG` [14], `SPDC` [37] have relied on this regularization step.[17]

With this new regularizer, one can show that the dual objective $D'(y)$ is $L_i \stackrel{\text{def}}{=} \frac{1}{n} + \frac{1}{\lambda_2 n^2}\|a_i\|_2^2$ smooth with respect to the $i$-th coordinate, as well as $1/n$-strongly convex with respect to the Euclidean norm. Therefore, we can again apply `RCDM`, `ACDM`, `APCG`, and our `NU_ACDM` to this objective $D'(y)$ and compare their performances. This is shown in Figure 3.

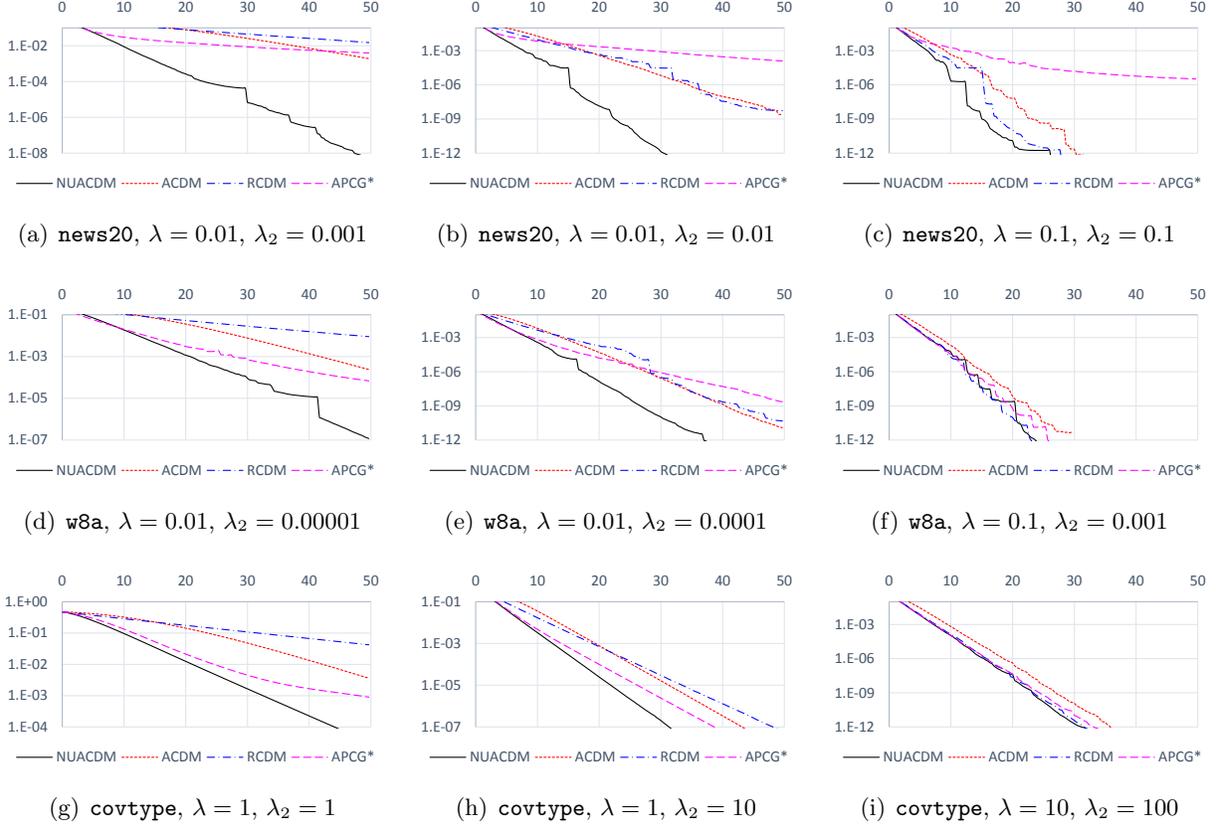

Figure 3: Performance Comparison for Lasso. The $y$ axis represents the objective distance to minimum, and the $x$ axis represents the number of passes to the dataset.

## 7.2 Experiments on Non-Strongly Convex Objectives

Consider problem (2.1) where $r(w) = \frac{\lambda}{2}\|w\|^2$ is the $\ell_2$ regularizer and each $\phi_i(\cdot)$ is some *non-smooth* loss function. In this case, the dual objective (2.2) becomes

$$\min_{y \in \mathbb{R}^n} \left\{ D(y) \stackrel{\text{def}}{=} \frac{1}{n}\sum_{i=1}^n \phi_i^*(y_i) + \frac{1}{2\lambda n^2}\Big\|\sum_{i=1}^n y_i a_i\Big\|_2^2 \right\} \ . \tag{7.2}$$

This $D(y)$ is not necessarily strongly convex because the penalty functions $\phi_i(\cdot)$ is not smooth. In this section, we conduct an experiment for the case when $\phi_i(\alpha) \stackrel{\text{def}}{=} \frac{1}{2}(\alpha - l_i)^2 + |\alpha - l_i|$ is an $\ell_2 - \ell_1$ penalty function. We call this ERM problem the $\ell_2 - \ell_1$ Penalty Regression.

As before, we know that $D(y)$ is $L_i \stackrel{\text{def}}{=} \frac{1}{n} + \frac{1}{\lambda n^2}\|a_i\|_2^2$ smooth with respect to the $i$-th coordinate, so we can apply `ACDM`, `RCDM`, `APCG` and our `NU_ACDM`<sup>ns</sup> directly to minimize $D(y)$. We choose $\beta = 0$ for `ACDM`, `RCDM`, and `NU_ACDM`<sup>ns</sup> in our experiment, and have to choose $\beta = 1$ for `APCG`. Our results are shown in Figure 4.

From these experiments, we see that again the theoretical speed-up factors in Table 2 are validated in practice. `NU_ACDM`<sup>ns</sup> has a clear advantage over its close relatives `ACDM` and `RCDM` when the coordinate smoothness parameters $L_i$ are very non-uniform (such as dataset `news20`), and when $\lambda$ is relatively small. In contrast, in the least non-uniform datasets covtype, ACDM performs even slightly faster than NU_ACDM$^{\mathtt{ns}}$. This is not surprising, because the theoretical speed-up in this case is only less than 0.3%.

---

[16]This reduction introduces a logarithmic factor loss in the running time, and was recently improved by [2].

[17]After this paper has been submitted for publication, a new accelerated method `Katyusha` was discovered to avoid this regularization step [1].

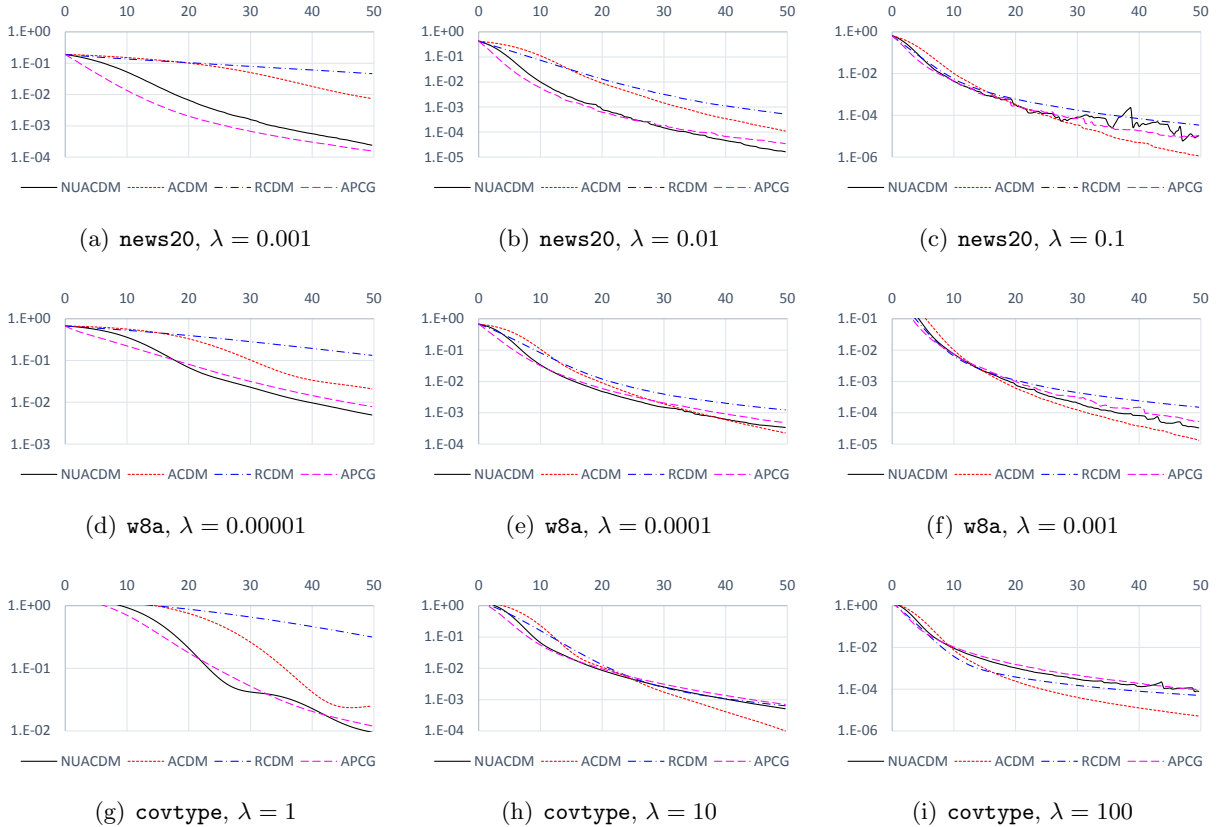

Figure 4: Performance Comparison for $\ell_1$-$\ell_2$ Penalty Regression. The $y$ axis represents the training objective distance to minimum, and the $x$ axis represents the number of passes to the dataset.

In contrast to the previous subsection, APCG (which uses $\beta = 1$) performs extremely well and similar to NU_ACDM$^{\mathtt{ns}}$ (which uses $\beta = 0$) in Figure 4. As we shall see in the next section, by taking a closer look at different choices of $\beta$ for non-strongly convex objectives, APCG is in fact analogous to the $\beta = 1$ case of NU_ACDM$^{\mathtt{ns}}$, but is slightly worse than NU_ACDM$^{\mathtt{ns}}$ for $\beta$ being between 0 and 0.8 for all the three datasets we are considering in this paper.

### 7.3 Dependence on $\beta$

As discussed in Remark 1.1, when dealing with a strongly convex objective $f(\cdot)$, we usually work with accelerated coordinate descent methods for Euclidean norm rather than $\mathsf{L}_\beta$ norms. However, the choice becomes less obvious for non-strongly convex objectives.

For instance, in Table 1, by comparing $T = \sum_i \sqrt{L_i/\varepsilon} \cdot \|x_0 - x^*\|$ for $\beta = 0$ and $T = n/\sqrt{\varepsilon} \cdot \|x_0 - x^*\|_{\mathsf{L}_1}$ for $\beta = 1$, it is not immediately clear which one is more preferable to the other. If one works with a standard machine learning boundedness assumption $\|x_0 - x^*\| \leq \Theta$ for some constant $\Theta$, then the convergence for the $\beta = 1$ case reduces to $T = n/\sqrt{\varepsilon} \cdot \|x_0 - x^*\|_{\mathsf{L}_1} \leq n \max_i \sqrt{L_i}/\varepsilon \cdot \Theta$ which is slower than that of the $\beta = 0$ case. However, in general, the best choice of $\beta$ depends on how the coordinates of the vector $x_0 - x^*$ scale with parameters $L_i$.

Nevertheless, we can perform a comparison *in practice* between difference choices of $\beta$. Focusing on the $\ell_1 - \ell_2$ Penalty Regression dual objective (7.2), we plot the performance of `NU_ACDM`$^{\text{ns}}$ with different $\beta$. From Figure 5, we conclude that smaller values of $\beta$ are perhaps more preferred to larger ones in practice. Not surprisingly, the performance difference becomes less significant for dataset `covtype`, because it has more uniform smoothness parameters $L_i$ than the other two datasets. Finally, we have included `APCG` in Figure 5 as well, and it has very similar performance comparing to `NU_ACDM`$^{\text{ns}}$ for $\beta = 1$. This confirms our theoretical finding in Table 1.

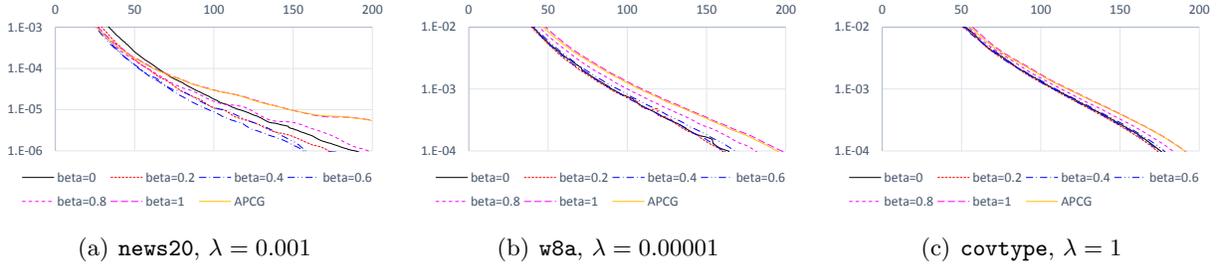

(a) `news20`, $\lambda = 0.001$  (b) `w8a`, $\lambda = 0.00001$  (c) `covtype`, $\lambda = 1$

Figure 5: Performance Comparison on `NU_ACDM`$^{\text{ns}}$ with $\beta = 0, 0.2, 0.4, 0.6, 0.8, 1$ and `APCG`.

|         | $r = 100\%$ | $r = 80\%$, | $r = 60\%$ | $r = 40\%$ | $r = 20\%$ | $r = 10\%$ |
|---------|-------------|-------------|------------|------------|------------|------------|
| Speed Up | 1           | 1.0992      | 1.2464     | 1.4025     | 1.6243     | 1.7379     |

Table 3: Theoretical Speed-Up Factors $\sqrt{n \sum_i L_i}/\left(\sum_i \sqrt{L_i}\right)$ of `NU_ACDM` over `ACDM` for linear systems $Ax = b$.

## 8  Experiments on Solving Linear Systems

We generate random linear systems $Ax = b$ and compare randomized Kaczmarz, `ACDM`, and `NU_ACDM`.

We choose $m = 300$ and $n = 100$, and generate each entry $A_{ij}$ uniformly at random in $[0, 1]$. We scale a fraction $r$ of $A$'s rows to have Euclidean norm 10, and the rest to have Euclidean norm 1. We generate a random vector $x$, compute $b = Ax$, and use each of the three algorithms to solve $x$ given $A$ and $b$.

Since the coordinate smoothness parameters depend on the Euclidean norm squares of $A$'s rows, we expect our `NU_ACDM` to have a greater speed up comparing to `ACDM` for small nonzeros values of $r$. We compute the theoretical speed up factors in Table 3.

In Figure 6, we see that both `NU_ACDM` and `ACDM` outperform the non-accelerated randomized Kaczmarz without surprise. Furthermore, `NU_ACDM` and `ACDM` are comparable for $r = 100\%$, and the out-performance indeed becomes more significant for smaller values of $r$.

## Acknowledgements

We thank Yin Tat Lee for helpful conversations and careful reading of a draft of this paper. ZA-Z is partially supported by a Microsoft Research Grant, no. 0518584. ZQ and PR would like to acknowledge support from the EPSRC Grant EP/K02325X/1, "Accelerated Coordinate Descent Methods for Big Data Optimization". PR also acknowledges support from the EPSRC Fellowship Grant EP/N005538/1, "Randomized Algorithms for Extreme Convex Optimization".

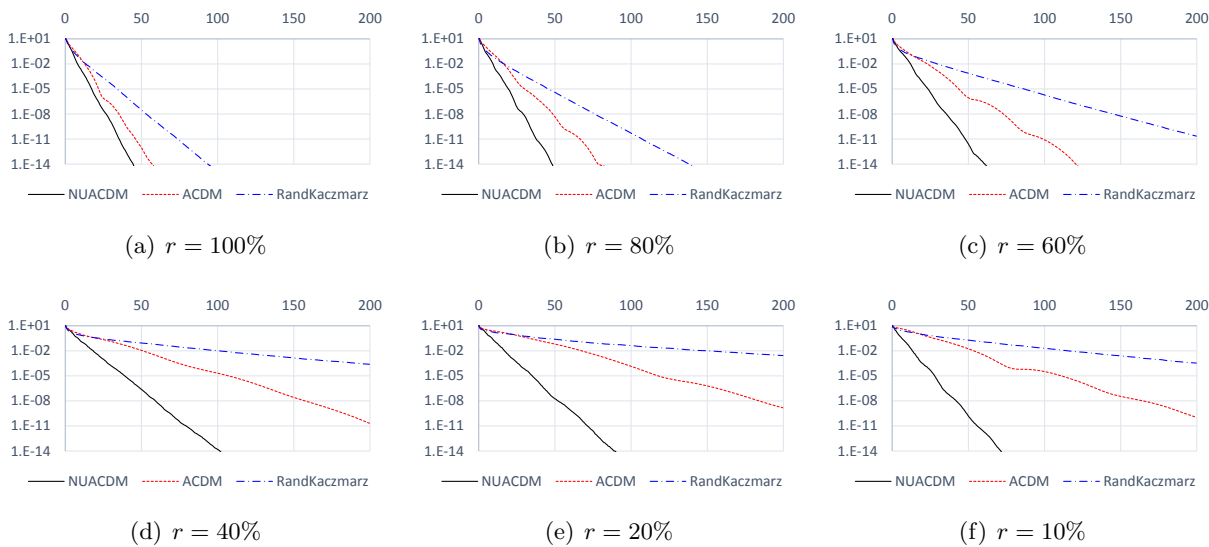

Figure 6: Performance Comparison on Solving $Ax = b$.

# APPENDIX

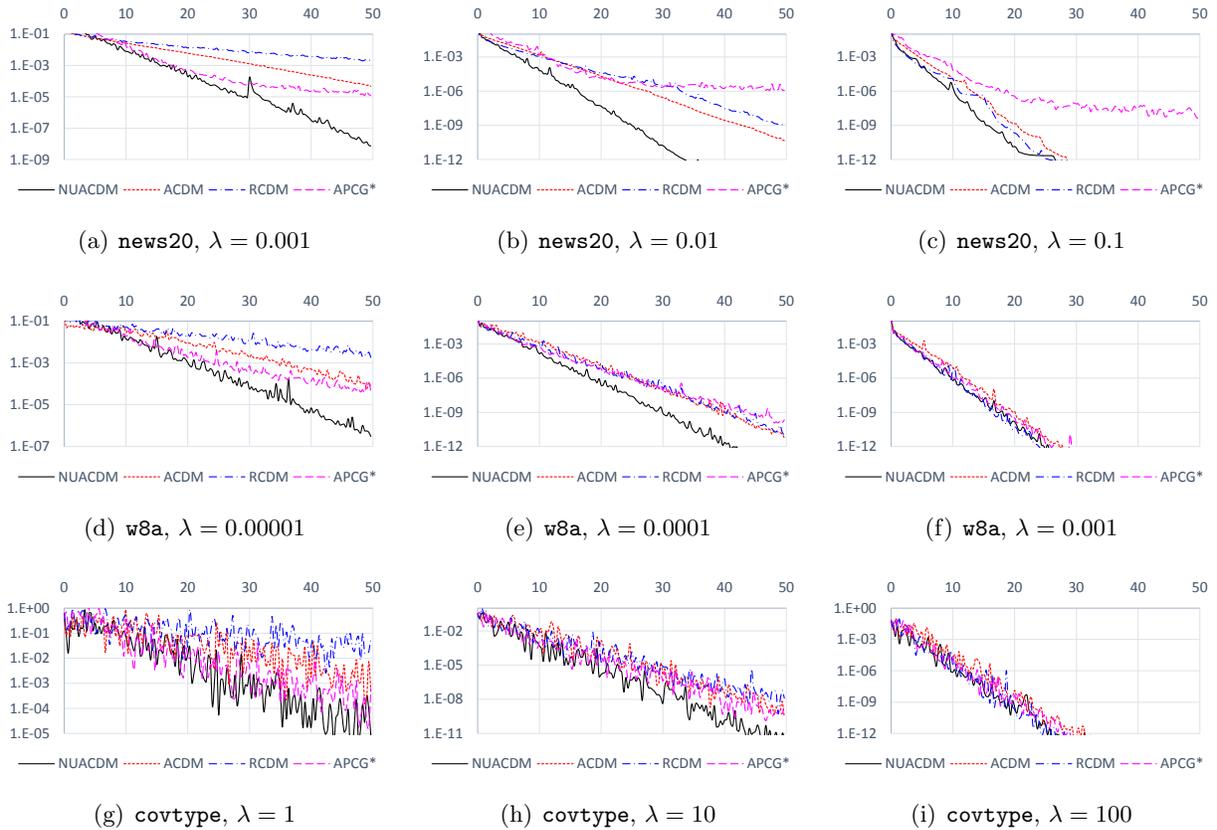

Figure 7: Performance Comparison for Ridge Regression (Primal). The $y$ axis represents the *primal* objective distance to minimum, and the $x$ axis represents the number of passes to the dataset.